\documentclass{article}
\usepackage{amssymb}
\usepackage{graphicx}
\usepackage{latexsym}
\def\part#1{\frac{\partial\phantom{q}}{\partial#1}}

\newenvironment{rmk}{\begin{trivlist}\item[]{\bf Remark:} }
{\end{trivlist}}

\newenvironment{prf}{\begin{trivlist}\item[]{\bf Proof:} }
{\hfill $\Box$ \end{trivlist}}

\newtheorem{thm}{Theorem}
\newtheorem{definition}{Definition}
\newtheorem{prp}[thm]{Proposition}

\newcommand{\lie}[1]{\mathfrak{#1}}
\newcommand\End{\mbox{End}}

\newcommand\Hom{\mbox{Hom}}

\newcommand\Imag{\mbox{Im}}

\newcommand\tr{\mbox{tr}}

\newcommand\ad{\mbox{ad}}

\newcommand\Real{\mbox{Re}}

\newcommand{\R}{\mathbf{R}}
\newcommand{\C}{\mathbf{C}}
\newcommand{\K}{\mathbf{H}}

\newcommand{\PP}{{\mathbf {\mathrm P}}}
\newcommand{\RP}{{\mathbf R}{\mathrm P}}
\newcommand{\OC}{\mathbf{O}}
\begin{document}

\title{$SL(2)$ over the  octonions}

\author{ Nigel Hitchin}

 \maketitle
 
 \begin{abstract} We interpret an open orbit in a $32$-dimensional representation space of $Spin(9,1)\times SL(2,\R)$ as a substitute for the non-existent group of invertible $2\times 2$ matrices over the octonions. The approach is via twistor geometry in eight dimensions.
  \end{abstract}

\section{Introduction}
The Lorentz groups $SO(n,1)$ have distinguished  properties in low dimensions, most notably through the special isomorphisms of their double covers. We have  $ Spin(2,1)\cong SL(2,\R), Spin(3,1)\cong SL(2,\C)$ and $Spin(5,1)\cong SL(2,\K)$ where $\K$ denotes the quaternions. In each case there is a corresponding geometrical statement: the group $SO(n,1)$ acts as conformal transformations of the sphere $S^{n-1}$ and then these isomorphisms are interpreted as M\"obius transformations of the projective line over the division algebras $\R,\C$ or $\K$. It is often said, especially in the physics literature (see \cite{CS},\cite{QF}), that the eight-sphere $S^8$ should be considered as  the projective line over the octonions $\OC$ and $ Spin(9,1)\cong SL(2,\OC).$ 
But there is something clearly wrong about this statement since the dimension of $SO(9,1)$ is $45$ and $2\times 2$ matrices with octonionic entries have only 32 dimensions. 

We propose here an alternative interpretation of $SL(2,\OC)$. We shall regard $GL(2,\OC)$ as an open orbit of the action on ${\mathbf S}\otimes \R^2$ of the group $Spin(9,1)\times GL(2,\R)$  where ${\mathbf S}$ is the 16-dimensional spin representation of $Spin(9,1)$. This is a natural object, one of the  irreducible prehomogeneous spaces in the list of Kimura and Sato \cite{Kim} and has the correct dimension. Furthermore, the stabilizer of a point is $G_2\times SL(2,\R)$ and the tangent space as a module for this group consists of  $2\times 2$ matrices with octonion entries. Then $SL(2,\OC)$ is the submanifold  given by setting a certain quartic polynomial   (the ``determinant") equal to $1$. It has an invariant inner product (the ``Killing form") and two natural subspaces which we interpret as $SU(2, \OC)$ (the ``maximal compact" one) and $SU(1,1, \OC)$.  There are then totally geodesic ``cosets" of $SL(2,\K)$ corresponding to quaternionic subalgebras of the octonions. 

The only problem is that these are not Lie groups. Nevertheless we show that they have many features in common with Lie groups.

Our chosen approach is via twistor theory in the Euclidean context of the conformal geometry of spheres. We begin with the familiar 4-dimensional story, pointing out some special issues about quaternionic  transformations, and then move to  the conformal geometry of the sphere $S^8$. 

The paper is based on a talk given in Oxford for Roger Penrose's 80th birthday on July 21st 2011.
\section{Twistors}
\subsection{Conformal geometry}
Here we recall the Euclidean version of Penrose's twistors. For the  four-dimensional version see \cite{AHS} for example. It allows us to define concrete models for the basic representations of the Lorentz group  in terms of the conformal differential geometry of the sphere. Recall that 
the quadric  $x_1^2+x_2^2+x_3^2+\cdots+x_{n+1}^2-x_0^2=0$ in projective space $\R \PP^{n+1}$ is diffeomorphic to the $n$-sphere since the homogeneous coordinate $x_0$ never vanishes, and its symmetry group  $SO(n+1,1)$ acts on $S^n$ as the group of orientation-preserving conformal transformations. 

Twistor theory concerns itself with spinors on $S^n$, or $\R^n$ using stereographic projection. When $n=2m$ is even there are  two spinor bundles $S^+$ and $S^-$ and two associated operators.

We have the {\emph{Dirac operator}} $D:S^{\pm}\rightarrow S^{\mp}$
$$D\psi=\sum_i e_i\cdot\nabla_i\psi$$
where $(e_1,\dots, e_n)$ is an orthonormal basis of tangent vectors which operate by Clifford multiplication  on the section $\psi$ of the spinor bundle. Recall that Clifford multiplication $\psi\mapsto u\cdot \psi$ by a tangent vector $u$ satisfies the basic identity $u\cdot v+v\cdot u=-2(u,v)1$.

There  is also the complementary \emph {twistor operator} $\bar D: S^{\pm}\rightarrow S^{\pm}\otimes T^*$ defined by 
$$\bar D\psi=\sum_i\nabla_i\psi\otimes e_i+\frac{1}{n} e_i\cdot D\psi\otimes e_i.$$
These expressions are written using the Levi-Civita connection of $\R^n$ with its flat metric $g$ but the operators make sense for any conformally equivalent metric -- if $g$ rescales by the function $\lambda^2$ (and so the \emph{conformal weight} of $g$ is 2) then $\psi$ scales by $\lambda^{(n-1)/2}$ for the Dirac operator and $\lambda^{-1/2}$ for the twistor operator.

The twistor equation $\bar D\psi=0$ is overdetermined but in flat space $\bar D:  S^{+}\rightarrow S^{+}\otimes T^*$ has a nullspace, the   space of \emph {twistors} ${\mathbf T}$.  For $\psi$ in $S^{-}$ this is the dual twistor space ${\mathbf T}^*$. A spinor field in ${\mathbf T}$ can be written on $\R^n$ as 
$$\psi={\mathbf x}\cdot \varphi^-+\varphi^+$$
where $\varphi^-$ and $\varphi^+$ are constant spinors and ${\mathbf x}\in \R^n$ is the position vector. Conformal invariance means that ${\mathbf T}$ is acted on by conformal transformations which lift to the spin structure. Then ${\mathbf T}$ and ${\mathbf T}^*$ are the two basic spin representations  of $Spin(n+1,1)$.

To make this  more precise we can consider the action of generators of the conformal group. If $F$ is a conformal  diffeomorphism its derivative $dF=\lambda P$ where $P$ is orthogonal. A twistor $\psi({\mathbf x})$ then transforms to $\lambda^{-1/2}\tilde P\psi(F({\mathbf x}))$ where $\tilde P$ is a lift of $P$ to the spin group. The action of a rotation is clear. For a\emph{ translation} by ${\mathbf a}$ we have 
$${\mathbf x}\cdot \varphi^-+\varphi^+\mapsto {\mathbf x}\cdot \varphi^-+({\mathbf a}\cdot \varphi^- +\varphi^+).$$
A \emph {reflection} in the hyperplane orthogonal to a unit vector ${\mathbf n}$ is orientation reversing and so maps $S^+$ to $S^-$ and vice-versa. The action is 
$${\mathbf x}\cdot \varphi^-+\varphi^+\mapsto {\mathbf n}\cdot (({\mathbf x-2({\mathbf x},{\mathbf n}){\mathbf n}})\cdot \varphi^- +\varphi^+)=-{\mathbf x}\cdot {\mathbf n}\cdot \varphi^-+{\mathbf n}\cdot \varphi^+$$
using the Clifford algebra identity.
An \emph{inversion} ${\mathbf x}\mapsto {\mathbf x}/r^2$ is also orientation reversing and the action is 
$${\mathbf x}\cdot \varphi^-+\varphi^+\mapsto r \frac{\mathbf x}{r}\cdot (\frac{1}{r^2}{\mathbf x}\cdot \varphi^-+\varphi^+)={\mathbf x}\cdot \varphi^+-\varphi^-.$$
 A \emph {dilation} ${\mathbf x}\mapsto \lambda {\mathbf x}$ gives 
 $${\mathbf x}\cdot \varphi^-+\varphi^+\mapsto \lambda^{-1/2}(\lambda{\mathbf x}\cdot \varphi^-+\varphi^+)={\mathbf x}\cdot \lambda^{1/2}\varphi^-+\lambda^{-1/2}\varphi^+.$$

There is another conformally invariant operator based on the trace-free Hessian of a function:
$$\bar\Delta f= \nabla^2 f-\frac{1}{n}(\Delta f)g$$
where $f$ rescales by $\lambda^{-1}$. The equation $\bar\Delta f=0$ has on $\R^n$ an $(n+2)$-dimensional space of solutions of the form 
\begin{equation}
f=a\Vert {\mathbf x}\Vert^2+({\mathbf x},{\mathbf b})+c
\label{scalars}
\end{equation}
Under the action of the conformal group $SO(n+1,1)$ this is the vector representation $V$, with the Lorentzian inner product 
$$L(f,f)=({\mathbf b},{\mathbf b})-4ac.$$

\begin{rmk}
If we rescale the flat metric $g$ by $f^{-2}$ where $\bar \Delta f=0$ we obtain a constant curvature metric on the complement of the zero locus of $f$.  If $L(f,f)=0$ then we have a Euclidean metric on the complement of the point, which we can regard as stereographic projection from that point on the sphere $S^n$. If $L(f,f)<0$ we have positive curvature. 
\end{rmk}

There is a relation between solutions to these  conformally invariant equations. The first defines the Clifford action of a Lorentzian vector in $V$ on the spinor space ${\mathbf T}$. If $\psi$ satisfies $\bar D\psi=0$ and $f$ satisfies $\bar\Delta f=0$ then
\begin{prp} \label{cliff} The linear map 
$$\psi\mapsto\frac{2}{n} f^{(n+2)/2}D(f^{-n/2}\psi)$$
is the Clifford action of $f\in V$ on $\psi\in {\mathbf T}$.
\end{prp}
\begin{prf} Note that the particular power of $f$ is chosen to give the right conformal weight to the Dirac operator. 

We set $\psi={\mathbf x}\cdot \varphi^-+\varphi^+$ and $f=a\Vert {\mathbf x}\Vert^2+({\mathbf x},{\mathbf b})+c$ and calculate:  
$$\frac{2}{n} f^{(n+2)/2}D(f^{-n/2}\psi)=-df\cdot \psi+\frac{2}{n} fD\psi=-(2a{\mathbf x}+{\mathbf b})\cdot \psi+\frac{2}{n}f\sum_ie_i\cdot e_i\cdot \varphi^-.$$
The right hand side is
$$-(2a{\mathbf x}+{\mathbf b})\cdot ({\mathbf x}\cdot \varphi^-+\varphi^+)- 2(a\Vert {\mathbf x}\Vert^2+({\mathbf x},{\mathbf b})+c)\varphi^-$$ and using the Clifford identities ${\mathbf x}\cdot {\mathbf x}=-\Vert {\mathbf x}\Vert^2$ and ${\mathbf x}\cdot {\mathbf b}+{\mathbf b}\cdot {\mathbf x}=-2({\mathbf b}, {\mathbf x})$ this gives 
$${\mathbf x}\cdot ({\mathbf b}\cdot\varphi^--2a\varphi^+)-({\mathbf b}\cdot \varphi^++2c\varphi^-)$$
which satisfies the twistor equation with the opposite chirality. 
Applying $f$ again  multiplies $\psi$ by $-L(f,f)$ and this gives the Clifford identity for the action of $\R^{n+1,1}$
\end{prf}

Depending on the dimension, the spinor bundles have a  Hermitian, or real,  inner product and  
$$\langle \psi,\psi\rangle =\langle {\mathbf x}\cdot\varphi^-+\varphi^+,{\mathbf x}\cdot\varphi^-+\varphi^+\rangle=\langle -{\mathbf x}\cdot{\mathbf x}\cdot\varphi^-,\varphi^-\rangle+2\Real \langle{\mathbf x}\cdot\varphi^-,\varphi^+\rangle+\langle\varphi^+,\varphi^+\rangle$$
since Clifford multiplication by ${\mathbf x}$ is skew-adjoint. But ${\mathbf x}\cdot{\mathbf x}=-\Vert {\mathbf x}\Vert^2$ then gives an expression of the form (\ref{scalars}) for $\langle \psi,\psi\rangle$, so another relation is that for two solutions $\psi_1,\psi_2$ of $\bar D \psi=0$, $f=\langle \psi_1,\psi_2\rangle$ satisfies $\bar\Delta f=0$.
\subsection{Dimensions $4$ and $8$}

On $\R^4$ the spin bundles $S^+$ and $S^-$ are trivial and  we shall use the same notation for  the representation spaces. These are the standard two-dimensional complex representations  of the two factors in $Spin(4)=SU(2)\times SU(2)$ and since $SU(2)$ is isomorphic to the unit quaternions they are one-dimensional quaternionic vector spaces. Writing $\psi={\mathbf x}\cdot \varphi^-+\varphi^+$ for constant spinors $\varphi^-, \varphi^+$ shows that  ${\mathbf T}$ is isomorphic  to pairs $(\varphi^-,\varphi^+)\in S^-\oplus S^+$ and is therefore a two-dimensional quaternionic vector space. The action of $Spin(5,1)$ preserves this quaternionic structure though not the direct sum  decomposition, as is evident from the action of translations.

There is another special property in four dimensions: a fixed nonzero spinor $\varphi^-$ defines an {\emph{isomorphism}} as real vector spaces between $\R^4$ and $S^+$ by the map ${\mathbf x}\mapsto {\mathbf x}\cdot \varphi^-$. It gives $\R^4$ a complex structure (this is the basis for the construction of the twistor space in \cite{AHS}). For us what this means is that ${\mathbf x}\cdot \varphi^-+\varphi^+= ({\mathbf x}-{\mathbf c})\cdot \varphi^-$ so that $\psi$ vanishes at a single point ${\mathbf x}={\mathbf c}$.  The function $f=\langle \psi,\psi\rangle$ satisfies $L(f,f)=0$ and in fact defines a flat metric on the complement of the point ${\mathbf x}={\mathbf c}$.

Since each  element $\psi$ in ${\mathbf T}$ defines a point ${\mathbf x}={\mathbf c}$ in $S^4=\R^4 \cup \{\infty\}$, this provides us with  a projection $p:{\mathbf T}\backslash \{0\}\rightarrow S^4$. The fibre is the one-dimensional quaternionic subspace spanned by 
 $(\varphi^-,-{\mathbf c}\cdot \varphi^- )\in S^-\oplus S^+$. We can regard this as the identification of $S^4$ with the quaternionic projective line $\K\PP^1$. Invariantly, the fibre over ${\mathbf c}$ is isomorphic to a fibre of the spinor bundle  $S^-$ over the sphere $S^4$.
 
In dimension $8$   the two spin bundles are real and  we have  a real inner product. We also have the phenomenon of triality \cite{Adams}.

The three $8$-dimensional spaces $\R^8,S^+,S^-$  have invariant maps 
$$\R^8\otimes S^+\rightarrow S^-,\quad \R^8\otimes S^-\rightarrow S^+,\quad S^+\otimes S^-\rightarrow \R^8.$$
The first two are Clifford multiplication, the last one is the adjoint of  the map $\R^8\rightarrow \Hom(S^+,S^-)\cong S^+\otimes S^-$ defined also by the Clifford action. We shall use the same notation $a\cdot b$ for each of these since they satisfy the same basic identities -- triality comes from an outer automorphism of $Spin(8)$ and permutes the roles of the three representation spaces.

The action of $Spin(9,1)$ on the space ${\mathbf T}$ of solutions $\psi={\mathbf x}\cdot \varphi^-+\varphi^+$ to $\bar D \psi=0$ is now a real $16$-dimensional representation. Triality interchanges the roles of vectors and spinors and so just as ${\mathbf x}$ gives an isomorphism from $S^+$ to $S^-$, so $\varphi^-$ gives an isomorphism from $\R^8$ to $S^+$ and hence an element of ${\mathbf T}$ can be written as $({\mathbf x}-{\mathbf c})\cdot \varphi^-$ as in the four-dimensional case  and defines a point in $S^8$. Moreover $f=(\psi,\psi)$ satisfies $L(f,f)=0$ and defines a flat metric on the complement of ${\mathbf c}$.

By analogy we should think of $S^8$ as the octonionic projective line but we are not yet in the position to make a convincing case for this. 
\section{$GL(2,\K)$ and $SL(2,\K)$}
\subsection{Quaternionic transformations}
A quaternionic vector space $U$ of dimension $m$ is invariantly  described as a complex vector space of dimension $2m$ together with an antilinear map $J:U\rightarrow U$ satisfying $J^2=-1$. 
 Then a complex linear transformation $A:U\rightarrow U$ is \emph{quaternionic} if $AJ=JA$ so that $A$ commutes with  $i,J$ and $K=iJ$. The quaternionic transformations are denoted by $\End_{\K}(U)$ and the invertible elements form the group $GL(U,\K)$.  Unlike complex linear transformations $\End_{\K}(U)$ is not a quaternionic vector space since for $q\in \K$, neither $qA$ nor $Aq$  commutes with $i$ and $J$. Nevertheless we can still represent $A$ by a matrix of quaternions.
 
 A quaternionic basis is a complex basis of the form 
 $\{e_1,Je_1,e_2,Je_2,\dots, e_m, Je_m\}$
and then if $A$ is quaternionic
$$A e_i=\sum_{j=1}^m B_{ji}e_j+C_{ji}Je_j\quad A Je_i=\sum_{j=1}^m-\bar C_{ji}e_j+ \bar B_{ji}Je_j$$
so $A$ is represented by an $m\times m$  quaternionic matrix $B+Cj$ acting on the left and commuting with multiplication on the right by $\K$.   

The vector space  $U\otimes_{\R}\R^m$ then has the quaternionic structure $J\otimes 1$ and is  an $m^2$-dimensional quaternionic vector space with a quaternionic action of $GL(U,\K)\times GL(m,\R)$. Choosing a basis it can be identified with $m\times m$ matrices with quaternionic entries. Hence the open set of invertible matrices is a free orbit of $GL(U,\K)$.  We can think of $U\otimes_{\R}\R^m$ as the $m$ columns of the quaternionic matrix.

\subsection{The determinant}\label{detsec}
Consider now the twistor space ${\mathbf T}$ for $\R^4$ as a $2$-dimensional quaternionic vector space and the space ${\mathbf T}\otimes_{\R}\R^2$. An element of this  vector space is a  pair 
$$\rho=(\psi_1,\psi_2)=({\mathbf x}\cdot \varphi_1^-+\varphi_1^+,{\mathbf x}\cdot \varphi_2^-+\varphi_2^+)$$ 
or a $2\times 2$ matrix where each entry lies in a one-dimensional quaternionic vector space
$$\rho=\pmatrix { \varphi_1^- & \varphi_1^+\cr
 \varphi_2^- & \varphi_2^+}$$
 
The Hermitian form $\langle \psi_a,\psi_b\rangle$ for $a,b=1,2$ defines  complex scalars $f_{ab}$ satisfying $\bar\Delta f_{ab}=0$ and 
then we can use the Lorentzian inner product to define 
$$\mu(\rho)=L(f_{11},f_{22})-L(f_{12},f_{21})$$
to obtain a real    homogeneous polynomial  of degree $4$.  Note that this is invariant under $Spin(5,1)\times SL(2,\R)$ and transforms by the factor $(\det P)^2$ for $P\in GL(2,\R)$.

If $A\in \End_{\K}(U)$ is a quaternionic transformation and $v$ is a complex eigenvector then $Av=\lambda v$ and 
$$AJv=JAv=J(\lambda v)=\bar\lambda Jv$$
so that $\bar\lambda$ is also an eigenvalue with eigenvector $Jv$. Hence the complex determinant of $A$ is real and non-negative. So the  determinant of a quaternionic $2\times 2$ matrix is a real polynomial of degree $4$.  
\begin{prp} \label{quatprop} If $\rho\in {\mathbf T}\otimes_{\R}\R^2$ is identified with a $2\times 2$ quaternionic transformation $A$, then $\mu(\rho)=-3\det A$.
\end{prp}
\begin{prf} Since $Spin(5,1)$ is simple its action via $GL(U,\K)$ preserves the determinant and by construction $\mu$ is invariant. 

Take $(\psi_1,\psi_2)\in {\mathbf T}\otimes_{\R}\R^2$. If one or other is zero then clearly both $\mu$ and the determinant are zero. If they are non-vanishing but their zeros are the same  then $\psi_1,\psi_2$ project to the same point in $S^4$, the quaternionic projective line, so $\psi_2=\psi_1 q$ for a quaternion $q=a+bj$ where $a,b$ are complex. Then $\psi_2=a\psi_1+b\psi_1j$ and so there is a complex linear relation among the 
columns of the $4\times 4$ complex matrix hence $\det A=0$.

On the other hand taking the common zero ${\mathbf c}$ to be $0$ we have  $\psi_1={\mathbf x}\cdot \varphi_1^-, \psi_2={\mathbf x}\cdot \varphi_2^-$ and the expressions $f_{12},f_{21},f_{11},f_{22}$ are all multiples of $r^2$. This gives  $\mu(\rho)=0$.

If the zeros are distinct they  can be taken to $0$ and $\infty$ in $S^4=\R^4\cup \{\infty\}$ by an element of $SO(5,1)$. 
This means $\psi_1={\mathbf x}\cdot \varphi^-, \psi_2=\varphi^+$. Then 
$$f_{11}=\langle {\mathbf x}\cdot \varphi^-,{\mathbf x}\cdot \varphi^-\rangle=r^2\Vert \varphi^-\Vert^2,\quad f_{22}=\Vert \varphi^+\Vert^2,\quad f_{12}=\langle {\mathbf x}\cdot \varphi^-,\varphi^+\rangle.$$
Using the identification of $S^+\otimes S^-$ with vectors, the latter term is the Euclidean inner product $({\mathbf x}, J\varphi^+\otimes \varphi^-)$ so now the Lorentzian inner products are given by 
$L(f_{11},f_{22})=-2\Vert \varphi^+\Vert^2\Vert \varphi^-\Vert^2$ and $$ L(f_{12},f_{21})=(J\varphi^+\otimes \varphi^-,J\varphi^+\otimes \varphi^-)=\Vert \varphi^+\Vert^2\Vert \varphi^-\Vert^2$$ 
and hence 
$$\mu(\rho)=L(f_{11},f_{22})-L(f_{12},f_{21})=-3\Vert \varphi^+\Vert^2\Vert \varphi^-\Vert^2$$
The corresponding quaternionic matrix $A$ is diagonal with entries $(q^+,q^-)=(\varphi^+,\varphi^-)$. As above a $1\times 1$ quaternion $q=a+bj$ is the complex matrix
$$\pmatrix{ a & b\cr -\bar b & \bar a}\,\, {\mathrm {with\,\, complex\,\, determinant }}\,\, a\bar a+b\bar b=q\bar q.$$
so $\det A=\Vert \varphi^+\Vert^2\Vert \varphi^-\Vert^2=-\mu(\rho)/3$.
\end{prf}

\begin{rmk} A corollary of the Proposition is the fact that $\det A$ defined as a complex determinant is independent of the particular complex structure in the quaternionic family $a_1i+a_2j+a_3k$ where $a_1^2+a_2^2+a_3^2=1$.
\end{rmk}
It follows from the above that we can identify $GL(2,\K)$ inside ${\mathbf T}\otimes \R^2$ as  the complement of the hypersurface $\mu(\rho)=0$ where $\mu$ is a $Spin(5,1)\times SL(2,\R)$-invariant quartic function. The group $SL(2,\K)$ is defined in this way by  $\mu(\rho)=-3$ and  Proposition \ref{quatprop} shows that this is identified with quaternionic transformations whose complex determinant is $1$. 

\section{$GL(2,{\mathbf O})$ and $SL(2,{\mathbf O})$} 
\subsection{The determinant}
By analogy with the quaternionic case we shall consider the space  ${\mathbf T}\otimes\R^2$ in the $8$-dimensional case under the action of $Spin(9,1)\times GL(2,\R)$ as a model for the space of $2\times 2$ matrices over the octonions. Before making this more specific we define   a real quartic invariant  which will play the role of the determinant. 

There are real inner products on the three spaces $\R^8, S^+,S^-$ which we will denote by the same symbol $\langle\,,\,\rangle $, and also write $\varphi^+\cdot \varphi^-\in \R^8$ for the invariant projection from $S^+\otimes S^-$ to $\R^8$, which under triality is equivalent to Clifford multiplication, which is a map from $\R^8\otimes S^{\pm}$ to $S^{\mp}$.

We take $(\psi_1,\psi_2)\in {\mathbf T}\otimes \R^2$ where $$\rho=(\psi_1,\psi_2)=({\mathbf x}\cdot \varphi_1^-+\varphi_1^+,{\mathbf x}\cdot \varphi_2^-+\varphi_2^+)$$ 
and define the quartic invariant as in Section \ref{detsec}, using the real inner products. 

Here 
$\langle\psi,\psi\rangle=\Vert {\mathbf x}\Vert^2\langle\varphi^-,\varphi^-\rangle+2\langle{\mathbf x}\cdot\varphi^-,\varphi^+\rangle+\langle \varphi^+,\varphi^+\rangle.$
Triality  defines a vector $\varphi^+\cdot\varphi^-$ 
which gives 
$$f_{ab}=\langle\psi_a,\psi_b\rangle=\Vert {\mathbf x}\Vert^2\langle\varphi_a^-,\varphi_b^-\rangle+\langle{\mathbf x},\varphi_a^+\cdot\varphi_b^-\rangle+\langle{\mathbf x},\varphi_b^+\cdot\varphi_a^-\rangle+\langle\varphi_a^+,\varphi_b^+\rangle$$
and then
\begin{eqnarray*}
L(f_{11},f_{22})
&=&4\langle\varphi_1^-\cdot\varphi_1^+,\varphi_2^-\cdot \varphi_2^+\rangle-2\Vert \varphi_1^-\Vert^2\Vert \varphi_2^+\Vert^2-2\Vert \varphi_2^-\Vert^2\Vert \varphi_1^+\Vert^2\\
L(f_{12},f_{21}) &=&  \Vert\varphi_1^+\cdot \varphi_2^-+\varphi_2^+\cdot \varphi_1^-\Vert^2 -4\langle\varphi_1^-,\varphi_2^-\rangle\langle\varphi_1^+,\varphi_2^+\rangle.
\end{eqnarray*}
But by skew-symmetry of the Clifford action $$\langle\varphi_1^-\cdot\varphi_2^+,\varphi_2^{-}\cdot \varphi_1^+\rangle=-\langle\varphi^-_2\cdot ( \varphi_1^-\cdot\varphi_2^+),\varphi_1^{+}\rangle$$
and the basic Clifford algebra identity $\phi\cdot\psi+\psi\cdot\phi=-2(\phi,\psi)1$ gives 
$$\langle\varphi^-_2\cdot (\varphi_1^-\cdot\varphi_2^+),\varphi_1^{+}\rangle=-\langle\varphi^-_1\cdot (\varphi_2^-\cdot\varphi_2^+),\varphi_1^{+}\rangle-2\langle\varphi_1^-,\varphi_2^-\rangle \langle\varphi_1^+,\varphi_2^+\rangle.$$
It follows that 
$$\langle\varphi_1^-\cdot\varphi_2^+,\varphi_2^{-}\cdot \varphi_1^+\rangle=-\langle\varphi_2^-\cdot\varphi_2^+,\varphi_1^-\cdot \varphi_1^{+}\rangle+2\langle\varphi_1^-,\varphi_2^-\rangle\langle\varphi_1^+,\varphi_2^+\rangle$$
which gives $L(f_{11},f_{22})=-2L(f_{12},f_{21})$. 

The quartic invariant is therefore given by 
$$\mu(\rho)=L(f_{11},f_{22})-L(f_{12},f_{21})=-3L(f_{12},f_{21}).$$

If we write 
$$\rho=\pmatrix{\varphi_1^- & \varphi_1^+\cr\varphi_2^- &\varphi_2^+}$$
and take a diagonal matrix $\varphi_1^+=0, \varphi_2^-=0$ then $\mu(\rho)=-3\Vert \varphi_2^+\Vert^2\Vert \varphi_1^-\Vert^2$ so we make the definition:
\begin{definition} The {\emph {determinant}} of $\rho\in {\mathbf T}\otimes \R^2$ is equal to $-\mu(\rho)/3$.
\end{definition}
\begin{rmk} The definition gives 
\begin{equation}
\det\rho=-\frac{1}{2}L(f_{11},f_{22})= \Vert \varphi_1^-\Vert^2\Vert \varphi_2^+\Vert^2+\Vert \varphi_2^-\Vert^2\Vert \varphi_1^+\Vert^2-2\langle\varphi_1^-\cdot\varphi_1^+,\varphi_2^-\cdot \varphi_2^+\rangle.
\label{detform}
\end{equation}
If we  evaluate the determinant in $4$ dimensions in the same way we obtain the known expression \cite{Park}, \cite{Wilk} 
$$\det A=\vert a\vert^2\vert d\vert^2+\vert b\vert^2\vert c\vert^2-2\Real[a\bar c d \bar b]$$
for the quaternionic matrix
$$A=\pmatrix{ a & b\cr
c & d}.$$
\end{rmk} 
Recall the projection $p:{\mathbf T}\backslash\{0\}\rightarrow S^8$. The following are properties  of the determinant: 
\begin{prp} \label{GL} Let $\rho=(\psi_1,\psi_2)$ be an element of  ${\mathbf T}\otimes \R^2$, then
\begin{enumerate}
\renewcommand{\theenumi}{\roman{enumi}}
\item
if $\psi_1$ and $\psi_2$ are nonzero then $\det \rho=0$ if and only if $p(\psi_1)=p(\psi_2)\in S^8$
\item
$Spin(9,1)\times GL(2,\R)$ acts transitively on the open set $\det(\rho)\ne 0$ in ${\mathbf T}\otimes \R^2$
\item
the stabilizer of $\rho$ where $\det \rho\ne 0$ is $G_2\times SL(2,\R)$.
\end{enumerate}
\end{prp}
\begin{rmk} The compact Lie group $G_2$ in the last item is the automorphism group of the octonions, so the proposition tells us that there is a natural identification of the tangent space of the subspace $\det \rho\ne 0$ with $2\times 2$ matrices with octonion entries. The first item justifies the interpretation of $S^8$ as the octonionic projective line. 
\end{rmk}
\begin{prf}

\noindent (i) As in Proposition \ref{quatprop} if $p(\psi_1)\ne p(\psi_2)$ we can take $\varphi^+_1=0,\varphi^-_2=0$ and from the formula  (\ref{detform})  this gives $\det \rho=\Vert \varphi_1^-\Vert^2\Vert \varphi_2^+\Vert^2\ne 0$. If $p(\psi_1)= p(\psi_2)=\infty\in S^8$ then 
$\psi_1= \varphi_1^+,\psi_2=  \varphi_2^+$ so that  $\varphi_1^-=0, \varphi_2^-=0$ and from the formula  we obtain $\det\rho=0$.

\noindent (ii) This is a known fact from \cite{Kim} but we give a proof using twistors. As before we can transform to $(\psi_1,\psi_2)=({\mathbf x}\cdot \varphi_1^-, \varphi_2^+)$.The two antipodal points are preserved by $SO(8)\subset SO(9,1)$. By  a diagonal $GL(2,\R)$-action we may assume  $\varphi_1^-, \varphi_2^+$ are of unit length. Then ${\mathbf e}=\varphi_1^-\cdot \varphi_2^+$ is a unit vector which can be rotated to $(1,0,0,\dots,0)$, stabilized by $SO(7)$.  Clifford multiplication by ${\mathbf e}$ identifies $S^+\cong S^-= S$ and then $\varphi_1^-, \varphi_2^+$ can be considered as spinors in the same space $S$, the spinor representation for $Spin(7)$.
But $Spin(7)$ acts transitively on the unit $7$-sphere in $S$. So ${\mathbf e}$ and $\varphi_1^-$ can be mapped to standard elements of $\R^8, S^-$ and then ${\mathbf e}=\varphi_1^-\cdot \varphi_2^+$ determines $\varphi_2^+$.

\noindent (iii) For $s\in \R$, translation by $s{\mathbf e}$ acts on  $\rho=({\mathbf x}\cdot \varphi^-,\varphi^+)$ as 
$$(\psi_1,\psi_2)\mapsto ({\mathbf x}\cdot \varphi^-+s{\mathbf e}\cdot \varphi^-,\varphi_+)=({\mathbf x}\cdot \varphi^-+s \varphi^+,\varphi^+)=(\psi_1+s\psi_2,\psi_2).$$
Dilation by $t$ acts as 
$$(\psi_1,\psi_2)\mapsto (t^{1/2}\psi_1,t^{-1/2}\psi_2)$$
and inversion followed by translation followed by inversion gives 
$$(\psi_1,\psi_2)\mapsto (\psi_1,\psi_2+u\psi_1).$$
These three types of conformal transformations generate a copy of $SL(2,\R)$ which commutes with the action of $Spin(7)$. But the action is also  realized  by the external  $SL(2,\R)$-action on $\rho$. Thus $\rho$ is preserved by  a diagonal copy of $SL(2,\R)$ in
$SL(2,\R)\times Spin(2,1)\times Spin(7)$. The remaining component of the stabilizer is $G_2\subset Spin(7)$, the stabilizer of a spinor, as above.
\end{prf}

The group $G_2$ is the group of automorphisms of the octonions. In fact, in the proof above, we had a triple  $(\varphi_2^+,\varphi_1^-,{\mathbf e})\in S^+ \times S^-\times  \R^8$ such that $({\mathbf e},\varphi_1^-\cdot \varphi_2^+)=1$. Then, as in  [\cite{Adams} theorem 15.14]  we can identify each eight-dimensional space with the octonions in such a way that each element of the triple $(x,y,z)=(\varphi_2^+,\varphi_1^-,{\mathbf e})$ maps to the identity $1\in \OC$ and 
\begin{equation}
\langle x,y\cdot z\rangle =\Real ((xy)z)
\label{triple}
\end{equation}

 From the Proposition  the action is transitive, so  the space  $\det \rho\ne 0$ is the homogeneous space  $Spin(9,1)\times GL(2,\R)/G_2\times SL(2,\R)$.  The Lie algebra  decomposes as 
$$\lie{so}(9,1) \oplus \lie{gl}(2)\cong (\lie{g}_2\oplus \lie {sl}(2))\oplus (\lie{gl}(2)\otimes \OC)$$
so the tangent space at a point is  $\lie{gl}(2)\otimes \OC$ as a $G_2\times SL(2,\R)$-module. This is naturally identified with the $32$-dimensional space of $2\times 2$ matrices over $\OC$. It seems natural then to say:

\begin{definition} The open orbit $\det\rho\ne 0$  is isomorphic to  $GL(2,\OC)$ and the submanifold $\det\rho=1$ is $SL(2,\OC)$.
\end{definition}

Of course, by the classification of Lie groups,  $GL(2,\OC)$ is not a group but nevertheless shares many features of Lie groups.

\subsection{Prehomogeneous spaces}
A {\emph {prehomogeneous space}} of a group $G$ is a representation space with an open orbit, so we have observed above an  example in  $Spin(9,1)\times GL(2,\R)$ acting on ${\mathbf T}\otimes \R^2$. This is one of the list in Kimura and Sato  \cite{Kim}, and other examples have been used as a basis for studying different geometries as in \cite{Hit1},\cite{Hit2},\cite{Hit3}. In these cases there is also a homogeneous function $\mu(\rho)$ which plays an important role in the geometry. For example when $GL(7,\R)$ acts on $\rho\in\Lambda^3 R^7$ then $\mu(\rho)$ is essentially the volume form of a $G_2$ structure.

Given that we are considering a prehomogeneous space we can strip away the twistor theory and just consider the linear algebra, but we lose this way the comparison with quaternions. So we consider ${\mathbf T}$ as simply the $16$-dimensional spin representation ${\mathbf S}$ of $Spin(9,1)$.  

Clifford multiplication by $v\in \R^{9,1}$ maps  ${\mathbf S}$ to the opposite spinor space, its dual ${\mathbf S}^*$. We have the canonical pairing  $\langle \psi,\psi'\rangle$ between $\psi\in {\mathbf S}$ and $\psi'\in {\mathbf S}^*$. Given $\psi_1,\psi_2\in {\mathbf S}$ we can therefore define
$\langle v\cdot\psi_1,\psi_2\rangle\in \R$, where $\langle \,,\,\rangle$ denotes the dual pairing, which is  {\emph symmetric} in $\psi_1,\psi_2$. Using the Lorentzian inner product $L$,  define  the symmetric bilinear expression $P(\psi_1,\psi_2)\in \R^{9,1}$ by
$$L(P(\psi_1,\psi_2),v)=\langle v\cdot\psi_1,\psi_2\rangle.$$
In particular we can define $Q(\psi)=P(\psi,\psi)$. This is a null vector, indeed ${\mathbf S}$ would have a quartic invariant $(Q(\psi), Q(\psi))$ otherwise. It defines a point in the quadric $S^8\subset \RP^9$ which we regard as the octonionic projective line.

The quartic  function $\mu$ on ${\mathbf S}\otimes \R^2$ is defined by
\begin{equation}
\mu(\rho)=L(Q(\psi_1),Q(\psi_2))
\label{deff}
\end{equation}
 where $\rho=(\psi_1,\psi_2)\in  {\mathbf S}\otimes\R^2$ (see also \cite{Gyo}). 
 
To see that there is an open orbit, note that under the action of $Spin(7)\times SL(2,\R)\subset Spin(9,1)$, ${\mathbf S}$ is expressed as a tensor product $S\otimes \R^2$ of spin representations of the two factors. Now  $G_2\subset Spin(7)$  is the stabilizer of a  spinor $\varphi$ in $S$ and $SL(2,\R)$ fixes a skew bilinear form $\epsilon \in \R^2\otimes\R^2$, thus $G_2\times SL(2,\R)$ stabilizes $\varphi \otimes \epsilon\in S\otimes\R^2\otimes\R^2={\mathbf S}\otimes \R^2$. It  has dimension $14+3=17$ and so its orbit under $Spin(9,1)\times GL(2,\R)$ has dimension $45+4-17=32=\dim {\mathbf S}\otimes \R^2$. 
 
 \subsection{The invariant metric}
 A simple Lie group has up to a scalar multiple a unique bi-invariant symmetric form on its Lie algebra, the tangent space at the identity. The Killing form $\tr (\ad A \,\ad B)$    is a canonical choice but for a linear Lie group where $A,B$ are matrices it is more useful to take $\tr AB.$ The determinant plays an essential role in defining the metric for a linear group. For $GL(n,\R)$ we have 
 $d\log \det P=\tr(P^{-1}dP)$
 and so, differentiating $P^{-1}$,  the Hessian of $\log \det P$ evaluated on $ A, B$ is 
 $$\tr(-P^{-1} AP^{-1} B)$$
 which is precisely $-\tr AB$ at the identity.

We therefore define a  metric on $SL(2,\OC)$ by the Hessian of   the logarithm of the function $\det \rho$. In fact since a scalar matrix $\lambda \in   GL(2,\R)$ takes $\log \rho$ to $\log\rho+2\log(\lambda^2)$ the Hessian is unchanged and so is invariantly defined on  $GL(2,\OC).$

\begin{prp} Let $A=A_0+A_1e_1+\dots+A_7e_7$ be a tangent vector to $SL(2,\OC)$ where $\{e_1,\dots,e_7\}$ is a standard basis for the imaginary octonions. Then the invariant inner product is 
$(A,A)=2(\sum_1^7\tr \,A_i^2-\tr\, A_0^2).$
It has signature $(22,9)$.
\end{prp}
\begin{prf}
For any function $f$ we have $$\nabla^2 \log f=\frac{\nabla^2 f}{f}-\frac{df\otimes df}{f^2}$$ and $df$ annihilates tangent vectors to $f=1$ so we only need the Hessian of $\det\rho$. We take this  at the point $\varphi^-_2=0,\varphi_1^+=0$ where $\varphi_1^-,\varphi_2^+$ are of unit length and 
 use the formula (\ref{detform})
$$\det\rho= \Vert \varphi_1^-\Vert^2\Vert \varphi_2^+\Vert^2+\Vert \varphi_2^-\Vert^2\Vert \varphi_1^+\Vert^2-2\langle \varphi_1^-\cdot\varphi_1^+,\varphi_2^-\cdot \varphi_2^+\rangle.$$
Then $\nabla \det \rho$ evaluated on a tangent vector is 
$$2[\langle \varphi_1^-,\dot \varphi_1^-\rangle+\langle \varphi_2^+,\dot  \varphi_2^+\rangle].$$
At the point under consideration $\varphi_1^-$ and $\varphi_2^+$ can be regarded, as  above,   as the identity in the octonions and the dotted terms as the entries in an octonionic matrix 
$$A=\pmatrix{ a & b\cr
c & d}$$
so $2[\langle \varphi_1^-,\dot \varphi_1^-\rangle+\langle\dot \varphi_2^+,\dot  \varphi_2^+\rangle]=2\Real(a+d)$ and the tangent space to $SL(2,\OC)$ consists of octonionic matrices whose trace is imaginary.

 The Hessian of $\det \rho$ is 
$$2[\langle\dot \varphi_1^-,\dot \varphi_1^-)+(\dot \varphi_2^+,\dot  \varphi_2^+\rangle+4\langle \varphi_1^-,\dot \varphi_1^-\rangle\langle \varphi_2^+,\dot  \varphi_2^+\rangle-2\langle\varphi_1^-\cdot \dot\varphi_1^+,\dot \varphi_2^-\cdot \varphi_2^+\rangle]$$
or, from (\ref{triple}), 
$2[a\bar a+d\bar d +4\Real(a)\Real(d)-2\Real(bc)]$. This is, using $a_0=-d_0$, 
$$\sum_{i=1}^72(a_i^2+d_i^2+2b_ic_i)- 2(2a_0^2+b_0c_0)=2(\sum_1^7\tr A_i^2-\tr A_0^2).$$
The signature is $7\times (3,1)+(1,2)=(22,9)$.
\end{prf}

\begin{rmk} Following the same procedure for the quaternions gives the metric 
$$2(\sum_1^3\tr A_i^2-\tr A_0^2)$$
whose signature is $3\times (3,1)+(1,2)=(10,5)$. In this case the quartic function is the genuine determinant which is invariant under left and right translation in $SL(2,\K)$ and so its Hessian defines a bi-invariant metric on the group. Since the group is simple this is  a multiple of the Killing form for $SL(2,\K)$ .
\end{rmk}

The quaternionic and octonionic case are clearly linked and the inclusion $Spin(5,1)\subset Spin(9,1)$ gives the following  ``cosets" of $SL(2,\K)$ in $SL(2,\OC)$.
\begin{prp} Through each point of $SL(2,\OC)$ there exists an eight-dimensional family of totally geodesic copies of $SL(2,\K)$, parametrized by the quaternionic subalgebras of the octonions. 
\end{prp}
\begin{prf}

 As above, $SL(2,\OC)$ is a homogeneous space of $Spin(9,1)\times SL(2,\R)$ with stabilizer $G_2\times SL(2,\R)\subset  Spin(7)\times Spin(2,1)\times SL(2,\R)$, the $SL(2,\R)$ factor being diagonally embedded. Any two orthogonal imaginary octonions generate a quaternion subalgebra of $\OC$.  The group $G_2$ acts transitively on $S^6$ and its stabilizer $SU(3)$ transitively on $S^5$ with stabilizer $SU(2)$. Hence a choice $i,j\in \K\subset \OC$ is parametrized by $G_2/SU(2)$.  But the automorphisms of $\K$ are $SO(3)$, and so the $8$-dimensional space $G_2/SO(4)$ parametrizes the quaternion subalgebras. Such a choice defines an $SO(4)$-invariant  decomposition of the tangent space to $SL(2,\OC)$ as
$\lie{gl}(2)\otimes (\K \oplus \R^4)$ together with the condition that the trace is imaginary. The action of $SO(4)$ on $\K \oplus \R^4$ is $1\oplus \Lambda^2_+\R^4\oplus \R^4.$

Consider the Lorentzian vector space $\R^{5,1}=\Imag \,\K\oplus \R^{2,1}\subset \R^{9,1}$ and its spin group $Spin(5,1)\subset Spin(9,1)$ and consider its action on $SL(2,\OC)$. Its intersection with the stabilizer is $Spin(3)\cap G_2$ corresponding to $SO(3)\subset SO(7)$. But the spin representation of $Spin(7)$ restricted to $Spin(3)$ is of the form $S\otimes \C^4$ which has no invariant vector. Since $G_2$ is the stabilizer of a spinor then the intersection is the identity and we have a free orbit of $Spin(5,1)=SL(2,\K)$. Its tangent space is $\lie{gl}(2)\otimes \K$ with real trace zero. 

To prove that these are totally geodesic, it is enough to show that the Hessian metric on $G/K=Spin(9,1)\times SL(2,\R)/G_2 \times SL(2,\R)$ is induced  from a bi-invariant metric on $G$, for the tangent space $\lie{gl}(2)\otimes \K$ is horizontal with respect to the orthogonal splitting. As a $G_2\times SL(2,\R)$-module the tangent space has  three irreducible components. If $V$, the imaginary octonions, is the 7-dimensional representation of $G_2$ then these are  
$$V.1\oplus \lie{sl}(2,\R)\oplus (\lie{sl}(2,\R)\otimes V).$$
We therefore have a three-dimensional family of invariant metrics. The group $Spin(9,1)\times SL(2,\R)$ has only a two-dimensional family, corresponding to the two simple factors. However, a bi-invariant metric on $Spin(9,1)$ restricts to one on $Spin(5,1)$ and we have seen that this gives $2(\sum_1^3\tr A_i^2-\tr A_0^2)$. This fixes the choice and  so the invariant metric is indeed a natural quotient metric.

\end{prf}

\subsection{Duality}

In general, when one has a relatively invariant functional $\mu$ on an open orbit in a vector space there is a corresponding open orbit and functional on the dual space and a nonlinear map from one orbit to the other. Differentiating $\mu$ at $\rho$ gives 
$$d \mu(\dot\rho)=\langle \hat \rho,\dot\rho\rangle$$
for a well-defined $\hat \rho$ in the dual space. An example is the action of $GL(7,\R)$ on $3$-forms $\rho\in \Lambda^3\R^7$ \cite{Hit1}. When $\rho$ is in an open orbit it defines a $G_2$ structure, and in particular a metric,  which gives $\hat\rho=\ast\rho \in  \Lambda^4\R^7$.

If we take $\mu(P)=\det P$ for $GL(n,\R)$ and use $\tr AB$ to identify the spaces of matrices with its dual then
$d\mu(A)=(\det P)\tr(P^{-1}A)$ and so $\hat P= (\det P) P^{-1}$. The orthogonal group $O(n)$ is then defined (for $n>2$) by an equation $\hat P=P^T$. We proceed next in a similar vein.

In our case $\rho\in {\mathbf T}\otimes \R^2$ and $\hat\rho\in ({\mathbf T}\otimes \R^2)^*$.   
Let $v$ be a vector in $\R^{9,1}$ such that $L(v,v)\ne 0$, then Clifford multiplication defines an isomorphism  $v:{\mathbf T}\rightarrow {\mathbf T}^*$ and the skew form $\epsilon$ on $\R^2$ defines an isomorphism with its dual.  Consider $\rho\in {\mathbf T}\otimes \R^2$ such that $\det \rho\ne 0$ and suppose
$$\hat\rho=(v\otimes \epsilon)\rho.$$ 
The space of such $\rho$ is acted on by the stabilizer of $v$ in $Spin(9,1)$ which is $Spin(9)$ if $L(v,v)<0$ and $Spin(8,1)$ if $L(v,v)>0$. 
\begin{prp} \label{compact} The subspace  $\{\rho\in SL(2,\OC): \hat\rho=(v\otimes \epsilon)\rho\}$ is isomorphic to 
\begin{enumerate}
\renewcommand{\theenumi}{\roman{enumi}}
\item
 $Spin(9)/G_2$ if $L(v,v)<0$
 \item
 $Spin(8,1)/G_2$ if $L(v,v)>0$. 
 \end{enumerate}
 The invariant metric has signature $(22,0)$ in the first case and $(14,8)$ in the second.
\end{prp} 
\begin{prf}

\noindent If $L(v,v)<0$, $v$ defines a positive curvature metric on $S^8$, which is conformally equivalent to $g/(1+r^2)^2$. So take $v$ to be represented by $f=(r^2+1)$ as a solution to $\bar\Delta f=0$. From Proposition  \ref{cliff} we have 
$v \cdot \psi_a=  -2({\mathbf x}\cdot \varphi_a^++\varphi_a^-)$ for $a=1,2$. Since $Spin(9)$  acts transitively and isometrically on $S^8$ then $\rho$ can be taken to a point with $\varphi_1^+=0$ 
and then 
$$(v\otimes \epsilon)\rho= 2(-{\mathbf x}\cdot \varphi_2^+-\varphi_2^-,\varphi_1^-).$$
Taking the derivative of $\det \rho$ at a point with $\varphi_1^+=0$ gives
$$\hat\rho=2(-\Vert \varphi_1^-\Vert^2{\mathbf x}\cdot \varphi_2^+,{\mathbf x}\cdot\varphi_1^-(\varphi_2^-\cdot\varphi_2^+)+\Vert\varphi_2^+\Vert^2\varphi_1^-).$$
Setting them equal gives $\varphi_2^-=0$ and $\Vert\varphi_1^-\Vert^2=\Vert\varphi_2^+\Vert^2=1$. This is the standard reference point we used in Proposition \ref{GL} and just as in that proof we have a transitive action on the pair of unit spinors with stabilizer $G_2$.

Differentiating both sides of $\hat\rho=(v\otimes \epsilon)\rho$ we obtain the equations for the tangent space to this submanifold (the derivative of $\hat \rho$ is essentially the Hessian). We obtain $\dot\varphi_2^+=\dot\varphi_2^++2\langle \dot\varphi_1^-,\varphi_1^-\rangle\varphi_2^+$ using $\Vert \varphi_1^-\Vert^2=1$ and $\varphi_2^-=\varphi_1^+=0$
and hence $\langle \dot\varphi_1^-,\varphi_1^-\rangle=0$ together with 
$$\dot\varphi_1^+=-\varphi_1^-\cdot(\dot\varphi_2^-\cdot \varphi_2^+)=\dot\varphi_2^-\cdot(\varphi_1^-\cdot \varphi_2^+)-2\langle\dot\varphi_2^-,\varphi_1^-\rangle  \varphi_2^+$$
and two more similar relations. 
The tangent space to $GL(2,\OC)$ at this point can be identified with $2\times 2$ matrices with octonion entries as we saw in Proposition \ref{GL} and the above relations read $\Real  (a)=0=\Real (d)$ and $b=c-2 \Real (c)=-\bar c$.
This  identifies  the $22$-dimensional tangent space as the subspace of  octonionic matrices 
$$A=\pmatrix{ a & b\cr
c & d}\qquad \bar A=-A^T.$$
Then the Hessian is 
$$2[a\bar a+d\bar d +4\Real(a)\Real(d)-2\Real(bc)]=2[a\bar a+d\bar d +2\Real(b\bar b)]$$
which is positive definite and so has signature $(22,0)$.

The procedure is entirely similar for the second part, using $f=(1-r^2)$. This defines the hyperbolic metric on the unit ball in $\R^8$ and $Spin(8,1)$ takes any point to the origin. 
\end{prf}
\begin{definition} The subspace  $Spin(9)/G_2\subset SL(2,\OC)$ we denote by $SU(2,\OC)$ and $Spin(8,1)/G_2 $ by $SU(1,1,\OC)$.
\end{definition}
The analogues in $SL(2,\K)$ are the maximal compact subgroup of  quaternionic unitary matrices $Sp(2)\cong Spin(5)$, and the indefinite version $Sp(1,1)\cong Spin(4,1)$ which are  genuine subgroups. 
\section{The compact space $SU(2,\OC)$}
The $22$-dimensional compact space $SU(2,\OC)=Spin(9)/G_2$ is often described as the octonionic Stiefel manifold $V_2(\OC)$. 
The proof of Proposition \ref{compact} gives a description of it without mention of the octonions. We showed there that we could normalize the pair of twistors to $(\psi_1,\psi_2)=({\mathbf x}\cdot \varphi^-,\varphi^+)$ which from the point of view of the sphere $S^8$ is a  point ${\mathbf x}=0$ with a unit spinor in $(S^+)_{\mathbf x}$  and one in $(S^-)_{\mathbf x}$. In other words it is the fibre product of the sphere bundles of $S^+$ and $S^-$ over $S^8$:
$$S^7\times S^7\rightarrow SU(2,\OC)\rightarrow S^8.$$ 

The action of $Spin(9)\subset Spin(9,1)$ on the unit sphere  in ${\mathbf S}$, with the metric given by the reduction to $Spin(9)$, is transitive with stabilizer $Spin(7)$. This  describes  another  fibration
\begin{equation}
S^7\rightarrow SU(2,\OC)\rightarrow S^{15}.
\label{15fibre}
\end{equation}
There are similar fibrations 
$S^3\times S^3\rightarrow Sp(2)\rightarrow S^4, S^3\rightarrow Sp(2)\rightarrow S^7$ in the quaternionic case.

The analogy with a compact Lie group is supported by the following:
\begin{prp} The space $SU(2,\OC)$ has the following properties:
\begin{enumerate}
\renewcommand{\theenumi}{\roman{enumi}}
\item
$SU(2,\OC)$ is a retraction of $SL(2,\OC)$,  
\item
its cohomology ring  is an exterior algebra on two generators,
\item
its tangent bundle is trivial.
\end{enumerate}
\end{prp}
\begin{prf}

\noindent (i) The metric $g/(1+r^2)^2$ on $\R^8$ completes to a constant positive curvature metric on the sphere $S^8$ and  the group $Spin(9,1)$ takes any two distinct points into antipodal points on $S^8$. This can be done in a systematic fashion: the length of the geodesic joining ${\mathbf a}$ to  $-{\mathbf a}$ is 
$$\int_{-a}^a\frac{dx}{(1+x^2)}= 2 \arctan a$$
where $a=\Vert {\mathbf a} \Vert$ so points are antipodal if $a=\pi/2$. The conformal map ${\mathbf x}\mapsto \lambda {\mathbf x}$ therefore takes the two points to antipodal ones if $\lambda=\pi/2a$. Then $\lambda$ is a smooth  function of ${\mathbf a}$ and is equal to $1$ if the points are antipodal. The action of a twistor is
$${\mathbf x}\cdot \varphi_-+\varphi_+\mapsto \lambda^{-1/2}(\lambda {\mathbf x}\cdot \varphi_-+\varphi_+)$$
or $(\varphi_-,\varphi_+)\mapsto (\lambda^{1/2}\varphi_-,\lambda^{-1/2}\varphi_+)$. This is an action by $SL(2,\R)$ which takes the point to the standard form for $SU(2,\OC)$.

\noindent (ii)  It follows from the structure of the fibration (\ref{15fibre}) that the cohomology is generated by classes in degree $7$ and $15$ whose product is the top degree class.  As pointed out in \cite{Kot}, this means that the product of two harmonic forms is harmonic, as in a Lie group. The cohomology ring has the same structure as that of  a rank $2$ simple Lie group, although for $Sp(2)$ the generators are in degrees $3,7$ and $SU(3)$ in degrees $3,5$.

\noindent (iii)  From \cite{Suth} to prove that $SU(2,\OC)$ has trivial tangent bundle $T$ it suffices, using the fact that it is a sphere bundle over a sphere, namely $S^{15}$,  to show that $1\oplus T$ is trivial. From the description of the tangent space above as $2\times 2$ octonionic matrices 
$$T\cong 1\oplus V \oplus V\oplus V$$
where $V$ is the rank $7$ bundle corresponding to the imaginary octonions. 

Considering the fibration  $p:SU(2,\OC)\rightarrow S^8$, a point in $SU(2,\OC)$ defines unit spinors $\varphi^+, \varphi^-$ at a point $x\in S^8$ and the unit vector $\varphi^+\cdot \varphi^-$ gives an isomorphism
$p^*TS^8\cong 1\oplus V$. For the tangent bundle of a sphere $TS^8\oplus 1$ is trivial, so, if $n$ denotes the trivial rank $n$ bundle we have $2\oplus V\cong 9$ and hence 
$$ 1\oplus (1\oplus V \oplus V\oplus V)\cong 9\oplus V\oplus V\cong (2\oplus V)\oplus (2\oplus V)\oplus 5\cong 23.$$ Hence $1\oplus T$ is trivial and then so is $T$ thanks to  [\cite{Suth}, theorem 1.3].
\end{prf}

\begin{rmk}
 The space $SU(1,1,\OC)$ fibres over the hyperbolic ball $B^8$ with fibre $S^7\times S^7$ and so retracts to $S^7\times S^7$. Compare this to $Sp(1,1)$ retracting to its maximal compact subgroup $Sp(1)\times Sp(1)\cong S^3\times S^3$.
\end{rmk}
\begin{rmk}
 The most obvious reason that $SU(2,\OC)$ is not a Lie group is that there is no degree $3$ cohomology class corresponding to the invariant three-form $B([X,Y]Z)$. In fact the exact homotopy sequence for the  fibration  over $S^{15}$ shows that $\pi_i(SU(2,\OC))=0$ for $1\le i \le 6$. Any compact simple Lie group $G$ has a copy of $SU(2)$ or $SO(3)$ which gives a generator of $\pi_3(G)$ but here this is replaced by the $7$-sphere.  The $3$-sphere which generates $\pi_3(Sp(2))$ is homotopically trivial in $S^7\subset SU(2,\OC)$.
\end{rmk}
\begin{rmk}
As in the non-compact situation of $SL(2,\OC)$ above, for each quaternion subalgebra of $\OC$ there is a totally geodesic copy of $Sp(2)$ through each point. Under the projection to the base $S^8$ each  such $Sp(2)$ projects to a totally geodesic $4$-sphere in $S^8$ giving the fibration $S^3\times S^3\rightarrow Sp(2)\rightarrow S^4$.  As we observed, there is an $8$-dimensional family of quaternionic subalgebras of $\OC$, and so the family of $Sp(2)$ subspaces has dimension $8+22-\dim Sp(2)=20$. But this is the dimension of  the Grassmannian  $SO(9)/SO(4)\times SO(5)$, so every totally geodesic $S^4\subset S^8$ arises this way and can be regarded as a quaternionic projective line inside the octonionic line.

\end{rmk}

\vskip 1cm
 \centerline{{\textit{Mathematical Institute, Woodstock Road, Oxford OX2 6GG}}}

 \end{document}